\documentclass[a4paper,10pt]{article}
\usepackage[latin1]{inputenc}
\usepackage[T1]{fontenc}
\usepackage{amsmath}
\usepackage{enumerate}
\usepackage{amsfonts}
\usepackage{amstext}
\usepackage{amsthm}
\usepackage[all]{xy}
\usepackage{geometry}
\usepackage{srcltx}
\usepackage{graphics}
\usepackage{epsfig}
\usepackage{subfigure}
\usepackage{slashed}
\usepackage{color}
\usepackage{bbm}
\usepackage{amssymb}
\usepackage{sidecap}
\usepackage{srcltx}
\usepackage{amscd}
\usepackage{comment}
\usepackage{tikz}

\DeclareMathOperator{\im}{im}

\DeclareMathOperator{\sgn}{sgn}

\DeclareMathOperator{\ind}{ind}

\DeclareMathOperator{\GL}{GL}

\newcommand{\R}{{\mathbb R}}
\newcommand{\Z}{{\mathbb Z}}

\newtheorem{theorem}{Theorem}[section]

\newtheorem{corollary}[theorem]{Corollary}

\theoremstyle{definition}

\newtheorem*{hypothesis*}{Hypothesis}
\theoremstyle{remark}

\newcommand{\cref}[1]{Cor.~\ref{#1}}

\newcommand{\href}[1]{Hyp.~\ref{#1}}

\title{A Topologist's Proof of Krasnosel'skii's Bifurcation Theorem}

\author{Robert Skiba and Nils Waterstraat}

\begin{document}
\date{}
\maketitle


\begin{abstract}
This note offers a proof of Krasnosel'skii's bifurcation theorem that bypasses the traditional machinery of nonlinear analysis: neither the Leray-Schauder degree nor the Lyapunov-Schmidt reduction is used. The argument rests only on elementary operator theory and straightforward vector-bundle topology, providing a topologist's perspective on a classical result in nonlinear functional analysis.
\end{abstract}

\section{Introduction}
Krasnosel'skii's bifurcation theorem provides one of the earliest and most influential criteria for detecting bifurcation from the trivial branch in nonlinear operator equations. Consider a family of nonlinear equations of the form 
\begin{align}\label{Krasnoequ}
u=\lambda (Ku+g(u))
\end{align}
on a Banach space $E$, depending on the parameter $\lambda\in\mathbb{R}$, where $K$ is a compact linear operator. Assume that the nonlinear map $g$ is Fr\' echet differentiable, satisfies $g(0)=0$ and $Dg(0)=0$, and maps bounded subsets of $E$ into relatively compact ones. Krasnosel'skii's Theorem asserts that if the compact operator $K$ has the eigenvalue $\frac{1}{\lambda_0}$ for some $\lambda_0\neq 0$, and if this eigenvalue has odd algebraic multiplicity, then $\lambda_0$ is a bifurcation point of the trivial branch, i.e. there exist elements $u_n\neq 0$ in $E$ and parameters $\lambda_n$ in $\mathbb{R}$ such that $u_n=\lambda_n(K u_n+g(u_n))$, $n\in\mathbb{N}$, and $(\lambda_n,u_n)$ converging to $(\lambda_0,0)$ in $\mathbb{R}\times E$.\\
Since its introduction in the fifties \cite{KrasnoI}, Krasnosel'skii's bifurcation theorem has become a workhorse across nonlinear science, especially in the analysis of qualitative transitions in mechanics, population dynamics, and reaction-diffusion systems, where a change in a parameter triggers the sudden appearance of new solution branches. In structural mechanics, for example, the theorem provides an explanation for the onset of Euler buckling, where the trivial straight configuration loses stability as a load crosses a critical eigenvalue. Its appeal lies in the fact that bifurcation can be certified purely from spectral information of a linear operator and compactness, without requiring variational structure or delicate transversality conditions.\\
There are different ways of proving this bifurcation theorem. Krasnosel'skii's original argument relies on the Leray-Schauder degree, which is defined for nonlinear compact perturbations of the identity as in \eqref{Krasnoequ} (see, e.g., \cite{Krasno}, \cite{Brown}). The idea is the following: if $\lambda_0$ is not a bifurcation point, then there is a neighbourhood $\Omega$ of $0\in E$ such that $u-\lambda(Ku+g(u))\neq 0$ for all $u\in\partial\Omega$ and all $\lambda$ sufficiently close to $\lambda_0$. By the compactness property of $g$, the Leray-Schauder degree of the map $I_E-\lambda(K+g)$ on $\Omega$ is well defined and remains constant in $\lambda$. As $K$ is compact, and hence $\frac{1}{\lambda}$ is not an eigenvalue of $K$ for $\lambda\neq\lambda_0$ close to $\lambda_0$, this allows one to compute the Leray-Schauder degree explicitly. This ultimately leads to a contradiction, since the degree changes as $\lambda$ passes through $\lambda_0$ when $\frac{1}{\lambda_0}$ has odd algebraic multiplicity.\\
A variant of this proof starts with a Lyapunov-Schmidt reduction. Using that $I_E-\lambda K$ is a Fredholm operator of index $0$, one reduces \eqref{Krasnoequ} to a family of nonlinear equations on $\ker(I_E-\lambda_0K)$. Since this kernel is finite-dimensional, one may replace the Leray-Schauder degree by the more elementary Brouwer degree. If no bifurcation occurs although $\frac{1}{\lambda_0}$ is an eigenvalue of $K$ of odd algebraic multiplicity, the Brouwer degree likewise yields a contradiction (see, e.g., \cite{Ize}).\\
The proof techniques outlined above are textbook examples of topological methods in nonlinear analysis. The purpose of this note is to present a proof of Krasnosel'skii's bifurcation theorem that differs from the standard ones, providing what may be called a topologist's proof. No prior knowledge of nonlinear analysis is required for our argument. Rather, the proof is designed for readers acquainted with the spectral theory of compact linear operators and with basic concepts from algebraic topology, especially vector bundles and (singular) homology. We do not wish to suggest that our proof is simpler or more elegant than the classical ones. Our hope is merely that avoiding the Leray-Schauder degree, and likewise the Lyapunov-Schmidt reduction, will make the argument attractive to readers with the appropriate prerequisites.

%%%%%%%%%%%%%%%%%%%%%%%%%%%%%%%%%%%%%%%%%%%%%%%%%%%%%%%%%%%%%%%%%%%%%%%%%%%%%%%%%%%%%%%%%%%%%%%%%%%%%%%%%%%%%%%%%%%%%%%%%%%%

\section{The Bifurcation Theorem}
Let $E$ be a real Banach space and $K=\{K_\lambda\}_{\lambda\in\mathbb{R}}$ a norm-continuous family of compact linear operators in $E$. In addition, we let $g:\mathbb{R}\times E\rightarrow E$ be a continuous map such that $g(\lambda,0)=0$ for all $\lambda$, each map $g_\lambda:=g(\lambda,\cdot):E\rightarrow E$ is Fr\' echet differentiable with $Dg_\lambda(0)=0$, and the map $(\lambda,u)\mapsto Dg_\lambda(u)$ from $\mathbb{R}\times E$ into $\mathcal{L}(E)$ is continuous. We consider the family of nonlinear operators
\begin{align}\label{G}
G(\lambda,u)=u-K_\lambda u+g(\lambda,u),\quad \lambda\in \mathbb{R},
\end{align}
which satisfies $G(\lambda,0)=0$ for all $\lambda\in\mathbb{R}$, and call $\lambda_0$ a bifurcation point of \eqref{G} if in every neighbourhood of $(\lambda_0,0)\in\mathbb{R}\times E$, there is some $(\lambda,u)$ such that $u\neq 0$ and $G(\lambda,u)=0$.\\
Note that $D_uG(\lambda,0)=I_E-K_\lambda$ by our assumptions on $g$. In particular, the implicit function theorem implies that $\lambda_0$ can be a bifurcation point only if $1\in\sigma(K_{\lambda_0})$.

\begin{theorem}\label{thm-main}
Let $a,b\in\mathbb{R}$ with $a<b$, and assume that $1$ is not an eigenvalue of $K_a$ and $K_b$. Suppose that there exists $0<r<1$ such that $\sigma(K_\lambda)\cap\partial D(0;r)=\emptyset$ for all $\lambda\in[a,b]$, where $D(0;r)$ denotes the closed disc of radius $r$ centered at $0\in\mathbb{C}$. Let $m(\lambda)$ denote the sum of the algebraic multiplicities of all real eigenvalues of $K_\lambda$ that are greater than $1$. If $m(a)+m(b)$ is odd, then there exists a bifurcation point $\lambda^\ast\in(a,b)$ for \eqref{G}.
\end{theorem}

Recall that if $K$ is a compact linear operator on $E$, then $0$ is the only point in the spectrum $\sigma(K)$ of $K$ that is possibly not isolated. Hence there exists $0<r<1$ such that $\sigma(K)\cap \partial D(0;r)=\emptyset$, and the complement $\mathbb{C}\setminus D(0;r)$ contains at most finitely many eigenvalues of $K$. Moreover, by continuity of isolated parts of the spectrum under norm perturbations, the same disk $D(0;r)$ has these properties for every compact operator $K'$ sufficiently close to $K$ in norm.
On the other hand, for the family $K_\lambda=\lambda K$, $\lambda\in[0,1]$, where $K$ is a fixed compact operator with spectral radius greater than $1$, no number $0<r<1$ as required in Theorem \ref{thm-main} can exist. Nevertheless, such families can still be handled by subdividing the path into smaller segments and applying the following corollary.

\begin{corollary}\label{cor-Krasnoselskii}
Let $a,b\in\mathbb{R}$ with $a<b$, and assume that there exists a unique parameter $\lambda_0\in(a,b)$ such that $1$ is an eigenvalue of $K_{\lambda_0}$, whereas $1\notin\sigma(K_\lambda)$ for all other $\lambda\in[a,b]$. 
Let $m(\lambda)$ denote the sum of the algebraic multiplicities of all real eigenvalues of $K_\lambda$ that are greater than $1$. If  $m(a)+m(b)$ is odd, then $\lambda_0$ is a bifurcation point for \eqref{G}.
\end{corollary}
Note that if $K_\lambda=\lambda K$ for some compact linear operator $K$, then by the previous corollary every nonzero parameter  $\lambda_0$ for which $\frac{1}{\lambda_0}$ is an eigenvalue of $K$ of odd algebraic multiplicity is a bifurcation point. This recovers precisely the conclusion of Krasonosel'skii's Theorem discussed in the introduction.\\
Let us finally note that neither Theorem \ref{thm-main} nor Corollary \ref{cor-Krasnoselskii} requires the compactness of $g$ in contrast to the classical Leray-Schauder approach (cf. \cite[\S 15.6]{Zeidler}). By using a Lyapunov-Schmidt reduction, compactness of $g$ is likewise unnecessary; however, in this case one usually assumes, as in Corollary \ref{cor-Krasnoselskii}, that $\lambda_0$ is isolated among all parameters  $\lambda\in\mathbb{R}$ for which $1$ is an eigenvalue of $K_\lambda$.\\
Rabinowitz later showed in \cite{Rabinowitz}, using the Leray-Schauder degree, that for compact 
$g$ in Krasnosel'skii's theorem, every bifurcation point actually gives rise 
to a branch of solutions which either becomes unbounded in $\mathbb{R}\times E$ 
or returns to the trivial branch $\mathbb{R}\times\{0\}$ (see \cite[\S 17]{Brown}). In our setting we do not assume compactness, yet our proof based on elementary topology nevertheless yields the existence of a local branch at the bifurcation points as follows.

\begin{theorem}\label{thm-mainII}
Under the assumptions of Theorem \ref{thm-main}, define 

\[\mathcal S:=\left\{(\lambda,u)\in[a,b]\times E:\, G(\lambda,u)=0,\ u\neq0\right\}.\]
If $m(a)+m(b)$ is odd, then for the bifurcation point $\lambda^\ast\in(a,b)$ of \eqref{G} in Theorem \ref{thm-main} (or likewise in Corollary \ref{cor-Krasnoselskii}), there exist $\rho_0>0,$ and a connected component $\mathcal C$ of
$\overline{\mathcal S}$ such that $(\lambda^\ast,0)\in\mathcal C$ and
\begin{equation}
		\label{eq:component-spheres}
		\mathcal C\cap
		\left\{
		(\lambda,u)\in[a,b]\times E:
		\|u\|=\rho
		\right\}
		\neq\emptyset
\end{equation}
for every $0<\rho<\rho_0.$
\end{theorem}
 We emphasize that more general bifurcation results for Fredholm operators of index $0$ were obtained by Fitzpatrick and Pejsachowicz \cite{FiPejsachowiczIII}, though their theory relies on substantially deeper tools from nonlinear functional analysis and therefore lies beyond the scope of the present paper.

%%%%%%%%%%%%%%%%%%%%%%%%%%%%%%%%%%%%%%%%%%%%%%%%%%%%%%%%%%%%%%%%%%%%%%%%%%%%%%%%%%%%%%%%%%%%%%%%%%%%%%%%%%%%%%%%%%%%%%%%%%%%%%%%%%%%%%%%%%%%%%%%%%%%%%%%%%%%%%%%%%%%%%%%%%%%%%%%%%%%%%%%%%%%%%%%%%%%%%%%%%%%%%%%%%%%%%%%%%%%%%%%%%%%%%%%%%%%%%%%%%%%%%%%%%%%%%%%%%%%%%%%%%%%%%%%%%%%%%%%%%%%%%%%%%%%%%%%%%%%%%%%%%%%%%%%%%%%%%%%%%%%%%%

\section{Proofs}
Henceforth we shorten notation by setting $L_\lambda=I_E-K_\lambda$, $\lambda\in\mathbb{R}$, which are Fredholm operators of index $0$ on $E$.

\subsection{Proof of Theorem \ref{thm-main}}
\subsubsection*{Step 1: A little operator theory - the bundles $\mathcal{V}$ and $\mathcal{W}$}
We begin by recalling some elementary facts from operator theory. For a compact linear operator $K$ on $E$, we define its spectrum by $\sigma(K):=\sigma(K^\mathbb{C})\subset\mathbb{C}$, where $K^\mathbb{C}$ denotes the complexification of $K$. Note that $\sigma(K)$ is invariant under complex conjugation. Let $D(0;r,R)\subset\mathbb{C}$ be an annulus centered at $0$ with inner radius $0<r<1$ and outer radius $R>1$, and assume that $\sigma(K)\cap \partial D(0;r,R)=\emptyset$. Then $D(0;r,R)$ contains only finitely many eigenvalues of $K$, and the associated spectral projection
\[P_{\partial D(0;r,R)}(K^\mathbb{C})=\frac{1}{2\pi i}\int_{\partial D(0;r,R)}(\mu I_{E^\mathbb{C}}-K^\mathbb{C})^{-1}\,d\mu\]
is of finite rank (see \cite[Thm. II.3.1]{Gohberg}). Its range $V^\mathbb{C}$ and kernel $W^\mathbb{C}$ are invariant subspaces of $K^\mathbb{C}$ and $\sigma(K^\mathbb{C}|_{V^\mathbb{C}})=\sigma(K)\cap D(0;r,R)$ (see \cite[Thm. I.2.2]{Gohberg}). Since the contour $\partial D(0;r,R)$ is invariant under complex conjugation, the projection $P_{\partial D(0;r,R)}(K^\mathbb{C})$ commutes with the canonical conjugation on $E^\mathbb{C}$. Hence there exists a unique real-linear projection $P_{r,R}(K)\in\mathcal{L}(E)$ such that $P_{r,R}(K)^\mathbb{C}=P_{\partial D(0;r,R)}(K^\mathbb{C})$.\\ 
We now consider the family $K=\{K_\lambda\}_{\lambda\in\mathbb{R}}$ and choose $R>1$ such that $R>\max\{\|K_\lambda\||\lambda\in[a,b]\}$. Since every spectral value $\mu\in\sigma(K_\lambda)$ satisfies $|\mu|<R$, our assumption on $r$ implies that $\sigma(K_\lambda)\cap\partial D(0;r,R)=\emptyset$ for all $\lambda\in[a,b]$. Hence, by \cite[Thm. II.4.1]{Gohberg}, the map $\lambda\rightarrow P_{r,R}(K_\lambda)$ is continuous from $[a,b]$ into $\mathcal{L}(E)$. As these projections have finite rank, they define a finite dimensional subbundle
\[\mathcal{V}=\{(\lambda,u)\in [a,b]\times E:\, u\in\im( P_{r,R}(K_\lambda))\}\]
of the product bundle $[a,b]\times E$ (see, e.g., \cite[\S 1]{Dupont}, \cite{indbundleIch}). By the preceding discussion, each fibre $\mathcal{V}_\lambda$ is invariant under $K_\lambda$, and therefore under $L_\lambda$ as well. Henceforth we shorten notation by setting $P_\lambda:=P_{r,R}(K_\lambda)$, $\lambda\in[a,b]$.\\
There is a complementary bundle $\mathcal{W}$, given fibrewise by the kernels of the projections $P_\lambda$, or equivalently by the ranges of the complementary projections $Q_\lambda:=I_E-P_\lambda$, which depend continuously on $\lambda$. Consequently, $\mathcal{V}\oplus\mathcal{W}=[a,b]\times E$ and each fibre $\mathcal{W}_\lambda$ is invariant under $L_\lambda$ as well.

\subsubsection*{Step 2: The map $F$ and finite-dimensional reduction via bundles}

 Since $L_a$ and $L_b$ are invertible and $\GL(E)$ is open, there exists $\varepsilon>0$ such that $L_\lambda$ is invertible for all $\lambda\in[a,a+2\varepsilon)\cup(b-2\varepsilon,b]$. Choose a continuous function $\chi:[a,b]\rightarrow[0,1]$ satisfying $\chi\equiv 0$ on $[a,a+\varepsilon)\cup(b-\varepsilon,b]$ and $\chi\equiv 1$ outside of $[a,a+2\varepsilon)\cup(b-2\varepsilon,b]$. Define $F:[a,b]\times E\rightarrow E$ by
\begin{align}\label{F}
F(\lambda,u)=(1-\chi(\lambda))L_{\lambda}u+G(\lambda,\chi(\lambda)\, u).
\end{align} 
This family has the same bifurcation points as \eqref{G}. Indeed, for every $\lambda$, $D_0F_{\lambda}=D_0G_{\lambda}=L_{\lambda}$, so $F$ has no bifurcation points in the region where $L_\lambda$ is invertible, i.e. on $[a,a+2\varepsilon)\cup(b-2\varepsilon,b]$. Outside this region we have $F_{\lambda}=G_{\lambda}$ and therefore the bifurcation points of $F$ coincide with those of $G$. The map $F$ has the additional property that $F_{\lambda}=L_{\lambda}$ for $\lambda=a,b$, a fact that will be needed at the end of the proof. From now on we work with the family $F:[a,b]\times E\rightarrow E$ and prove Theorem \ref{thm-main} by showing that, under the given assumptions, $F$ possesses a bifurcation point.\\
The following finite-dimensional reduction is inspired by an argument of Pejsachowicz \cite{JacoboTMNAI}. Define
\[F^1:[a,b]\times E\rightarrow\mathcal{V},\, F^1_\lambda=P_\lambda F_\lambda,\qquad F^2:[a,b]\times E\rightarrow\mathcal{W},\, F^2_\lambda=Q_\lambda F_\lambda,\]
and consider
\[\Phi:[a,b]\times E\rightarrow E,\quad\Phi(\lambda,u)=F^2_\lambda(u)+P_\lambda u.\]
For each $\lambda$, the map $\Phi_\lambda:E\rightarrow E$ is differentiable, and $(D\Phi_\lambda)(u)$ depends continuously on $(\lambda,u)\in [a,b]\times E$. Moreover,
\begin{align}\label{PQPhi}
(D_u\Phi)(\lambda,0)=D(Q_\lambda F_\lambda)(0)+P_{\lambda}=Q_\lambda L_{\lambda}+P_{\lambda}.
\end{align}
Since $\im(Q_{\lambda}L_{\lambda})=\mathcal{W}_{\lambda}$ and $\ker(Q_{\lambda}L_{\lambda})=\mathcal{V}_{\lambda}$, which are respectively the kernel and range of $P_{\lambda}$, the operator $(D_u\Phi)(\lambda,0)$ is invertible. By a simple modification of the inverse mapping theorem to include a parameter (see \cite[Lem. 2.1]{AmbrosettiProdi}), one obtains, for every $\lambda_0\in[a,b]$, an interval $I\subset[a,b]$ containing $\lambda_0$, a neighbourhood $\mathcal{U}$ of $I\times\{0\}$ in $I\times E$, and a neighbourhood $V$ of $0$ in $E$, such that
\[\Psi:\mathcal{U}\rightarrow I\times V,\quad \Psi(\lambda,u)=(\lambda,\Phi(\lambda,u))\]
is a homeomorphism. Since $[a,b]$ is compact, there exist $\delta>0$ and an open neighbourhood $\tilde{\mathcal{U}}$ of $[a,b]\times\{0\}$ in $[a,b]\times E$ such that $\Psi:\tilde{\mathcal{U}}\rightarrow[a,b]\times B_\delta(0)$ is a homeomorphism, where $B_\delta(0)$ denotes the open ball of radius $\delta$ in $E$.\\
Finally, it follows from the definition of $\Phi$ that $F^2(\Psi^{-1}(\lambda,u))=0$ for all $u\in\mathcal{V}_\lambda\cap B_\delta(0)$. Hence for $u\in\mathcal{V}_\lambda\cap B_\delta(0)$ we see that $F(\Psi^{-1}(\lambda,u))=0$ if and only if $F^1(\Psi^{-1}(\lambda,u))=0$. Consequently, we have reduced the bifurcation problem to the finite-dimensional vector bundle $\mathcal{V}$ over $[a,b]$.

%%%%%%%%%%%%%%%%%%%%%%%%%%%%%%%%%%%%%%%%%%%%%%%%%%%%%%%%%%%%%%%%%%%%%%%%%%%%%%%%%%%%%%%%%%%%%%%%%%%%%%%%%%%%%%%%%%%%%%%%%%%%

\subsubsection*{Step 3: A little topology}
Henceforth we let $k$ denote the dimension of the bundle $\mathcal{V}$, that is, the total number of eigenvalues lying in the annulus $D(0;r,R)$ for any $\lambda$. As $[a,b]$ is a contractible topological space, the bundle $\mathcal{V}$ admits a global trivialisation  $\varphi:\mathcal{V}\rightarrow [a,b]\times\mathbb{R}^k$ (see, e.g., \cite[Rem 1.51]{Dupont}). For sufficiently small $\rho>0$, we can define the map 
\begin{align}\label{f}
	f:=\varphi\circ (F^1\circ\Psi^{-1})\circ\varphi^{-1}:[a,b]\times B_\rho(0)\rightarrow [a,b]\times\mathbb{R}^k.
\end{align}
Let us now consider the endpoints $\lambda=a$ and $\lambda=b$. Then $F(\lambda,u)=L_\lambda u$, and thus $\Phi_\lambda(u)=F^2(\lambda,u)+P_\lambda u=Q_\lambda L_\lambda u+P_\lambda u$. As in \eqref{PQPhi}, this map is an isomorphism of $E$, acting as the identity on $\mathcal{V}_\lambda$. Consequently, $(F^1\circ\Psi^{-1})(\lambda,u)=(\lambda,L_\lambda u)$ for $\lambda\in\{a,b\}$ and all $u\in\mathcal{V}_\lambda\cap B_\delta(0)$. In particular, for $\lambda=a,b$, the fibre maps $f_\lambda$ induce group homomorphisms
\begin{align}\label{indhom}
(f_\lambda)_\ast:H_k(B_\rho(0),B_\rho(0)\setminus\{0\})\rightarrow H_k(\mathbb{R}^k,\mathbb{R}^k\setminus\{0\}),
\end{align}
where $(H_n)_{n\geq 0}$ denotes singular homology with integer coefficients. For the linear maps $f_a$ and $f_b$, the homomorphism $(f_\lambda)_\ast$ is multiplication by $\sgn(\det(f_\lambda))=\sgn(\det(\varphi_\lambda\circ L_\lambda\circ\varphi^{-1}_\lambda))=\sgn(\det(L_\lambda|_{\mathcal{V}_\lambda}))$ (cf., e.g., \cite[Prop. 10.5.3]{Dieck}). 
Finally,
\begin{align}\label{det}
\det(L_\lambda|_{\mathcal{V}_\lambda})=\prod^k_{i=1}(1-\mu_i(\lambda)),
\end{align}
where $\mu_1(\lambda),\ldots,\mu_k(\lambda)$ are the eigenvalues of the compact operator $K_\lambda$ in $D(0;r,R)$, counted with multiplicity. Non-real eigenvalues occur in complex conjugate pairs and therefore do not contribute to the sign of the determinant \eqref{det}, nor do eigenvalues lying in the intervals $(-R,-r)$ and $(r,1)$. Hence the sign of \eqref{det} for $\lambda=a,b$ depends only on the eigenvalues of $K_a$ and $K_b$ in $(1,R)$. Since $R$ exceeds the modulus of every eigenvalue of $K_a$ and $K_b$, it follows that 
\begin{align*}
\sgn(\det(L_\lambda|_{\mathcal{V}_\lambda}))=(-1)^{m(\lambda)},\quad \lambda=a,b.
\end{align*}
Consequently, the assumption that $m(a)+m(b)$ is odd implies that
\begin{align}\label{inequHom}
(f_a)_\ast\neq (f_b)_\ast.
\end{align}
To prove Theorem \ref{thm-main}, we now assume that there is no bifurcation point of $F$, aiming for a contradiction. By taking $\rho>0$ smaller in \eqref{f}, we may assume that the zero set of $f$ in $[a,b]\times B_\rho(0)$ is $[a,b]\times\{0\}$. Then for every $\lambda\in[a,b]$ we can regard $f_\lambda$ as a map $f_\lambda:(B_\rho(0),B_\rho(0)\setminus\{0\})\rightarrow (\mathbb{R}^k,\mathbb{R}^k\setminus\{0\})$, which yields $(f_a)_\ast=(f_b)_\ast$ by the homotopy invariance of singular homology and thus the desired contradiction by \eqref{inequHom}.

%%%%%%%%%%%%%%%%%%%%%%%%%%%%%%%%%%%%%%%%%%%%%%%%%%%%%%%%%%%%%%%%%%%%%%%%%%%%%%%%%%%%%%%%%%%%%%%%%%%%%%%%%%%%%%%%%%%%%%%%%%%%%%%%%%%%%%%%%%%%%%%%%%%%%%%%%%%%%%%%%%%%%%%%%%%%%%%%%%%%%%%%%%%%%%%%%%%%%%%%%%%%%%%%%%%%%%%%%%%%%%%%%%%%%%%%%%%%%%%%%%%%%%%%%%%%%%%%%%%%%%%%%%%%%%

\subsection{Proof of Theorem \ref{thm-mainII}}
In this section, the axioms of singular homology are no longer sufficient, and we therefore begin by recalling the homological definition of the Brouwer degree, where our main reference is \cite[\S 8-9]{Brown}. If $U\subset\R^k$ is bounded and open, and $\emptyset\neq\Sigma\subset U$ is compact, then there is a canonical embedding $j:(U,U\setminus\Sigma)\rightarrow (S^k,S^k\setminus\Sigma)$ where we consider $\mathbb{R}^k$ as subset of $S^k$ by stereographic projection. Note that $j$ induces isomorphisms in homology by the excision property. Composing with the monomorphism $i_\ast:H_k(S^k)\rightarrow H_k(S^k,S^k\setminus\Sigma)$ coming from the canonical inclusion $i$ gives a homomorphism $j^{-1}_\ast i_\ast:H_k(S^k)\rightarrow H_k(U,U\setminus\Sigma)$ and thus a non-trivial element $\nu(U,\Sigma)\in H_k(U,U\setminus\Sigma)$ as image of the canonical generator of $H_k(S^k)$ under this map.\\   
If now $f:\overline U\to\R^k$ is continuous, and $0\notin f(\partial U),$ then the compact zero set $\Sigma$ of $f$ in $U$ determines an integer $\ind(f,U)\in\Z$ by $f_\ast(\nu(U,\Sigma))=\ind(f,U)\,\nu_k$, where $f_\ast:H_k(U,U\setminus\Sigma)\rightarrow H_k(\R^k,\R^k\setminus\{0\})$ and $\nu_k$ is the canonical generator of $H_k(\R^k,\R^k\setminus\{0\})$. Note that if $0\in U$ and $0$ is the only zero of $f$ in $U$, then
$\ind(f,U)$ is the unique integer such that $f_*(\mu_U)=\ind(f,U)\,\mu,$
where $f_*:H_k(U,U\setminus\{0\})\rightarrow H_k(\R^k,\R^k\setminus\{0\})$
is the induced homomorphism and $\mu_U$ and $\mu$ are the canonical generators
determined by the standard orientation of $\R^k$. This is precisely the
interpretation used above in the proof of Theorem~\ref{thm-main}. In particular, if $f$ is linear and invertible, then $\ind(f,U)=\sgn\det f.$ Finally, if $\Omega\subset[a,b]\times\R^k$ is a bounded relatively open set and
$f:\overline\Omega\to\R^k$ continuous such that $\Omega_\lambda:=\{v\in\R^k:(\lambda,v)\in\Omega\}\neq\emptyset$ for all $\lambda\in [a,b]$ as well as  $f(\lambda,v)\neq0$ for every $(\lambda,v)\in\partial\Omega$, then 
\begin{align}\label{HomInvDeg}
	\ind(f_a,\Omega_a)
	=
	\ind(f_b,\Omega_b).
\end{align}
To prove Theorem \ref{thm-mainII} we now return to the proof of Theorem \ref{thm-main} and consider the map $f$ introduced in \eqref{f}. As before, we let $k\in\mathbb{N}$ be the dimension of the bundle $\mathcal{V}$ over $[a,b]$. Henceforth, we assume that $\rho>0$ is so small that $f$ is defined on $[a,b]\times\overline B_\rho(0)$ and has no non-trivial zeros with $\|u\|\leq\rho$ in the endpoint regions $[a,a+2\varepsilon)\cup(b-2\varepsilon,b],$ where we use that $F$ has no bifurcation point in $[a,a+2\varepsilon)\cup(b-2\varepsilon,b]$ by Step 2 in the proof of Theorem \ref{thm-main}. Note that with this choice of $\rho$ the endpoint homomorphisms in \eqref{indhom} are represented on $B_\rho(0)$.
Set $Q_\rho:=[a,b]\times \overline B_\rho(0)$ and define

\[
M_\rho:=
\overline{\left\{(\lambda,u)\in Q_\rho:f(\lambda,u)=0,\ u\neq0\right\}}.
\]
Since $Q_\rho$ is compact, $M_\rho$ is compact as well.  Set
\[
A_\rho:= M_\rho\cap \left([a,b]\times\{0\}\right)
\quad
\text{and}
\quad
\Gamma_\rho:= M_\rho\cap \left\{(\lambda,u): \|u\|=\rho \right\}.
\]
We first note that we have already shown in Step 3 in the proof of Theorem \ref{thm-main} that $A_\rho\neq\emptyset$. Now suppose that $\Gamma_\rho=\emptyset.$
Then
\[
f(\lambda,u)\neq0
\quad
\text{for every }
\lambda\in[a,b],
\quad
\|u\|=\rho,
\]
and \eqref{HomInvDeg} applied to $\Omega=[a,b]\times B_\rho(0)$ shows that the endpoint degrees coincide. 
Since at the endpoints the origin is the only zero, this yields $(f_a)_*=(f_b)_*,$ contrary to \eqref{inequHom}. Hence
$\Gamma_\rho\neq\emptyset.$\\
We now prove that some connected component of $M_\rho$ meets both
$A_\rho$ and $\Gamma_\rho$. Assume, to the contrary, that no connected component of $M_\rho$
meets both sets. By Whyburn's separation lemma (see, e.g., \cite[\S 14]{Brown}), there exist disjoint compact sets
$M_A,M_\Gamma\subset M_\rho$ such that
\begin{equation}
	\label{eq:Whyburn-decomposition}
	M_\rho=M_A\cup M_\Gamma,
	\quad
	A_\rho\subset M_A,
	\quad
	\Gamma_\rho\subset M_\Gamma.
\end{equation}
Let $T:=[a,b]\times\{0\}$ be the trivial branch and define $\widehat M_A:=M_A\cup T.$ Since $M_\rho\cap T=A_\rho\subset M_A,$
the compact sets $\widehat M_A$ and $M_\Gamma$ are disjoint. Furthermore, since every point of $M_\rho$ with norm $\rho$ belongs
to $\Gamma_\rho\subset M_\Gamma$, the set $\widehat M_A$ is disjoint from $[a,b]\times\partial B_\rho(0).$ Hence we can choose a relatively open set
$\Omega\subset[a,b]\times B_\rho(0)$ such that $\widehat M_A\subset\Omega,$ $\overline\Omega\cap M_\Gamma=\emptyset.$ We claim that
\begin{equation}
	\label{eq:nozeros-boundary-Omega}
	f(\lambda,u)\neq0
	\quad
	\text{for every }(\lambda,u)\in\partial\Omega.
\end{equation}
Indeed, every trivial zero belongs to $T\subset\Omega$, whereas every nontrivial zero in $Q_\rho$ belongs to $M_\rho=M_A\cup M_\Gamma$.
The set $M_A$ is contained in $\Omega$ and $M_\Gamma$ is disjoint from $\overline\Omega$. Thus no zero can lie on $\partial\Omega$.
For $\Omega_\lambda:=\{z\in\R^k:(\lambda,z)\in\Omega\},$ \eqref{HomInvDeg} now gives
$\ind(f_a,\Omega_a)=\ind(f_b,\Omega_b).$ Since the origin is the only zero at $a$ and $b$, excision reduces
these indices to the homomorphisms $(f_a)_*$ and  $(f_b)_*,$ respectively. Consequently, $(f_a)_*=(f_b)_*,$ contrary to \eqref{inequHom}.
Thus our assumption was false, and hence there exists a connected component $K_\rho$ of $M_\rho$ such that
\begin{equation}
	\label{eq:Krho-connects}
	K_\rho\cap A_\rho\neq\emptyset,
	\quad
	K_\rho\cap\Gamma_\rho\neq\emptyset.
\end{equation}
Since $M_\rho$ is compact, $K_\rho$ is compact as well. Thus $K_\rho$ indeed is a continuum.\\
We now consider the map $\Psi^{-1}\circ\varphi^{-1}:[a,b]\times \overline{B}_\rho(0)\rightarrow [a,b]\times E$, which is a homeomorphism onto its range, and define $\widetilde{K}_\rho=(\Psi^{-1}\circ\varphi^{-1})(K_{\rho})$. Clearly, as $K_{\rho}$ is compact and connected, this also holds for $\widetilde{K}_\rho$. Moreover, as in Step 2 above, the definitions of $\Psi$, $F^1$ and $F^2$ imply that $F(\lambda,u)=0$ for every $(\lambda,u)\in\widetilde{K}_\rho$. By the choice of $\rho$, all non-trivial points of $K_\rho$ lie outside the endpoint regions $[a,a+2\varepsilon)\cup(b-2\varepsilon,b]$ on which $F$ and $G$ coincide. Consequently, $\widetilde{K}_\rho\subset\overline{\mathcal{S}}$. As $K_\rho\cap A_\rho\neq\emptyset$, we see that $\widetilde{K}_\rho$ contains a bifurcation point $(\lambda^\ast,0)$ of $G$ for some $\lambda^\ast\in[a+2\varepsilon,b-2\varepsilon]$. Finally, as $K_\rho\cap\Gamma_\rho\neq\emptyset$, $\widetilde{K}_\rho$ contains a point $(\lambda_1,u_1)$ with $u_1\neq0.$ Set $\rho_0:=\|u_1\|>0.$ Consider the continuous map
$q:\widetilde K_\rho\to\R,$ $q(\lambda,u):=\|u\|.$ Since $\widetilde K_\rho$ is connected, its image $q(\widetilde K_\rho)$ is a connected subset of $\R$. Moreover,
$0=q(\lambda^\ast,0)$ and $\rho_0=q(\lambda_1,u_1).$ Therefore $[0,\rho_0]\subset q(\widetilde K_\rho).$ Consequently, for every
$0<\rho'<\rho_0$ there exists $(\lambda,u)\in\widetilde K_\rho$ such that $\|u\|=\rho'.$ Finally, let $\mathcal C$ be the connected component of
$\overline{\mathcal S}$ which contains $\widetilde K_\rho$. Then $(\lambda^\ast,0)\in\mathcal C,$ and, since
$\widetilde K_\rho\subset\mathcal C,$ we have
\[
\mathcal C\cap
\left\{
(\lambda,u):
\|u\|=\rho'
\right\}
\neq\emptyset
\]
for every $0<\rho'<\rho_0.$ Thus \eqref{eq:component-spheres} holds, which completes the proof.

%%%%%%%%%%%%%%%%%%%%%%%%%%%%%%%%%%%%%%%%%%%%%%%%%%%%%%%%%%%%%%%%%%%%%%%%%%%%%%%%%%%%%%%%%%%%%%%%%%%%%%%%%%%%%%%%%%%%%%%%%%%%%%%%%%%%%%%%%%%%%%%%%%%%%%%%%%%%%%%%%%%%%%%%%%%%%%%%%%%%%%%%%%%%%%%%%%%%%%%%%%%%%%%%%%%%%%%%%%%%%%%%%%%%%%%%%%%%%%%%%%%%%%%%%%%%%%%%%%%%%%%%%%%%%%

\subsection{Proof of Corollary \ref{cor-Krasnoselskii}}
We first show that, if $\alpha,\beta\in[a,b]$ are such that $L_\lambda$ is invertible for all $\lambda\in[\alpha,\beta]$, then $(-1)^{m(\alpha)}=(-1)^{m(\beta)}$. To prove this, fix $\lambda_1\in[\alpha,\beta]$ and let $0<r<1$ be such that $\partial D(0;r)\cap\sigma(K_{\lambda_1})=\emptyset$. By continuity of isolated parts of spectra \cite[Thm. II.4.1]{Gohberg}, there exists $\varepsilon>0$ such that $\partial D(0;r)\cap\sigma(K_\lambda)=\emptyset$ for all $\lambda\in[\lambda_1-\varepsilon,\lambda_1+\varepsilon]$. Let $R>1$ be larger than the modulus of every eigenvalue of the operators $K_\lambda$ for this interval and consider the bundle $\mathcal{V}$ constructed for the family $\{L_\lambda\}_{\lambda\in[\lambda_1-\varepsilon,\lambda_1+\varepsilon]}$ as in the first step of the proof of Theorem \ref{thm-main}. Then \eqref{det} holds for all $\lambda\in[\lambda_1-\varepsilon,\lambda_1+\varepsilon]$, and as above we obtain $(-1)^{m(\lambda_1-\varepsilon)}=(-1)^{m(\lambda_1+\varepsilon)} $. Since $[\alpha,\beta]$ is connected, this shows that $(-1)^{m(\alpha)}=(-1)^{m(\beta)}$.\\
Now let $\lambda_0\in(a,b)$ be as in Corollary \ref{cor-Krasnoselskii}, and choose again $r>0$ and $\varepsilon>0$ such that $\partial D(0;r)\cap\sigma(K_{\lambda})=\emptyset$ for all $\lambda\in[\lambda_0-\varepsilon,\lambda_0+\varepsilon]$. By the previous paragraph, $(-1)^{m(a)}=(-1)^{m(\lambda_0-\varepsilon)}$ and $(-1)^{m(b)}=(-1)^{m(\lambda_0+\varepsilon)} $. Consequently, $m(\lambda_0-\varepsilon)+m(\lambda_0+\varepsilon)$ is odd, and the claim follows from Theorem \ref{thm-main}.

%%%%%%%%%%%%%%%%%%%%%%%%%%%%%%%%%%%%%%%%%%%%%%%%%%%%%%%%%%%%%%%%%%%%%%%%%%%%%%%%%%%%%%%%%%%%%%%%%%%%%%%%%%%%%%%%%%%%%%%%%%%%%%%%%%%%%%%%%%%%%%%%%%%%%%%%%%%%%%%%%%%%%%%%%%%%%%%%%%%%%%%%%%%%%%%%

\vspace{1cm}
Robert Skiba\\
Faculty of Mathematics and Computer Science\\
Nicolaus Copernicus University in Toru\'n\\
Chopina 12/18\\
87-100 Torun\\
Poland\\
E-mail: robert.skiba@mat.umk.pl

\vspace{1cm}
Nils Waterstraat\\
Martin-Luther-Universit\"at Halle-Wittenberg\\
Naturwissenschaftliche Fakult\"at II\\
Institut f\"ur Mathematik\\
06099 Halle (Saale)\\
Germany\\
E-mail: nils.waterstraat@mathematik.uni-halle.de

\end{document}